\documentclass{article}
\usepackage{enumitem}
\usepackage{mathtools}
\usepackage{amsthm}
\usepackage[margin=2.65cm]{geometry}
\usepackage{tikz}
\usetikzlibrary{decorations.markings}
\usepackage[width=0.875\textwidth,font=small]{caption}
\usepackage{titling}

\begin{document}

\title{\huge Supersaturation in union-closed families of sets}
\author{\Large Christopher Bouchard}
\date{}
\maketitle

\vspace{-0.75cm}

\abstract{\vspace{-0.05cm} \noindent Let $k$ and $n$ be positive integers such that $2 \leq k \leq n+1$. We prove that the number of $k$-chains in a union-closed family with universe $[n]$ and size $m$ is minimized when its member sets are largest possible. We also show that, whenever the minimum is nonzero and $m>n$, there are no other minimizing families.}

\vspace{0.875cm}

\section*{1. Introduction}

Let $\mathcal{A}$ be a finite family of distinct finite sets, at least one of which is nonempty, and denote by $[n]=\{1,\cdots,n\}$ its \textit{universe} $U(\mathcal{A}) \coloneqq \bigcup_{A \in \mathcal{A}}A$. We refer to the family $\binom{[n]}{p} \coloneqq \{X \subseteq [n] \colon |X| = p\}$ as a \textit{layer} of the power set $2^{[n]} \coloneqq \{X \colon X \subseteq [n]\}$ for any integer $p$ such that $0 \leq p \leq n$. A \textit{chain} $\mathcal{C}$ in $\mathcal{A}$ is a subfamily of $\mathcal{A}$ such that $X_1 , X_2 \in \mathcal{C}$ implies that $X_1 \subseteq X_2 \ \lor \ X_2 \subseteq X_1$.

Sperner's foundational theorem of \cite{12} states that the maximum size of any family $\mathcal{A}$ with no chain of size two is achieved uniquely by a middle layer of $2^{[n]}$. A well-known generalization, given by Erd\H os in \cite{6}, establishes for $2 \leq k \leq n+1$ the maximum size of any $\mathcal{A}$ without a $k$\textit{-chain} (chain of size $k$) to be the sum of the $k-1$ largest binomial coefficients of $n$, which is achieved only by a union of $k-1$ middle layers. Kleitman in \cite{9} conjectured that (for $m$ greater than Erd\H os' size bound) the minimum number of $k$-chains in a family $\mathcal{A}$ of size $m$ is attained by concentrating member sets toward the center of the power set, having verified this for $k=2$ (a case itself building on a conjecture of Erd\H os and Katona). Kleitman's supersaturation conjecture for $2 \leq k \leq n+1$ was subsequently studied in \cite{1}, \cite{4}, and \cite{5}, with full proof and extremal characterization provided by Samotij in \cite{11}.
 
On another note, a family $\mathcal{A}$ is \textit{union-closed} if $X_1,X_2 \in \mathcal{A}$ implies that $X_1 \cup X_2 \in \mathcal{A}$. Most research on such families has dealt with the union-closed sets conjecture, also referred to as Frankl's conjecture, which states that if $\mathcal{A}$ is union-closed, then there exists an element from $[n]$ that is in at least half of its member sets. A breakthrough determining existence of an element in a constant fraction of member sets can be found in \cite{8}. For an earlier survey of various other results, see \cite{3}. 

In \cite{2}, we established the union-closed version of Erd\H os' bound, i.e. that for $2 \leq k \leq n+1$, the maximum size of a union-closed family $\mathcal{A}$ with no chain of size $k$ is the sum of the first $k-1$ binomial coefficients of $n$, attained uniquely by the union of the top $k-1$ layers of the power set. In the present work, we prove the corresponding supersaturation result, namely that if $m$ exceeds the sum of the first $k-1$ binomial coefficients of $n$, then the minimum number of $k$-chains in a union-closed family $\mathcal{A}$ of size $m$ is achieved by setting member sets as close as possible to the top of the power set (thus as large as possible). Moreover, we show that if $m$ also exceeds $n$, then there are no other extremal families.

We refer to a family $\mathcal{A}$ with universe $[n]$ as \textit{top-aligned} if $|A_1| \geq |A_2|$ for all $A_1 \in \mathcal{A}$ and $A_2 \in 2^{[n]} \setminus \mathcal{A}$. The main theorem is stated as follows:

\vspace{0.375cm}

\noindent \textbf{Theorem.} \textit{For any positive integers} $k$ \textit{and} $n$ \textit{such that} $2 \leq k \leq n+1$\textit{, the number of} $k$\textit{-chains in a union-closed family with universe} $[n]$ \textit{and size} $m$ \textit{is minimized whenever the family is top-aligned. Further, if the minimum is nonzero and $m>n$, then there are no other minimizing families.}

\section*{2. Proof of the main theorem}

\noindent Denote by $\textsf{C}(\mathcal{F},k)$ the collection of $k$-chains in a family $\mathcal{F}$, and by $c(k,m,n)$ the minimum possible number of $k$-chains in a union-closed family with universe $[n]$ and size $m$, recalling that $2 \leq k \leq n+1$. Let $\mathcal{A}_0$ be any non-top-aligned union-closed family with universe $[n]$ and size $m$. To prove the main theorem, we will consecutively apply two operators $T_1$ and $T_2$ to $\mathcal{A}_0$ in order to obtain the transformed family $\mathcal{A}_{\min}=(T_2 \circ T_1)(\mathcal{A}_0)$. The operators $T_1$ and $T_2$ are chosen such that:

\medskip

\begin{enumerate}[label=\roman*., align=left]
\item $|\mathcal{A}_{\min}|=|\mathcal{A}_0| = m$;
\item $U(\mathcal{A}_{\min}) = U(\mathcal{A}_0) = [n]$; 
\item $\mathcal{A}_{\min}$ is top-aligned;
\item $|\textsf{C}(\mathcal{A}_{\min},k)| \leq |\textsf{C}(\mathcal{A}_0,k)|$ (and also $|\textsf{C}(\mathcal{A}_{\min},k)| < |\textsf{C}(\mathcal{A}_0,k)|$ if $m>n$ and $c(k,m,n) > 0$).
\end{enumerate}

\medskip

\noindent Thus, any union-closed family that is not top-aligned can be transformed into a top-aligned family with same size, identical universe, and decreased number of $k$-chains (strictly decreased when family size exceeds universe size and minimum possible number of $k$-chains is nonzero), which proves the theorem.

\subsection*{2.1. Applying $T_1$ to $\mathcal{A}_0$}

Let $\mathcal{A}$ be a union-closed family of sets with universe $[n]$. We first consider, as described in \cite{3}, the \textit{up-compression} $u_x(\mathcal{A})$ of $\mathcal{A}$ with respect to an element $x \in [n]$, defined as

\vspace{0.1cm}

\[\hspace{-0.1cm}u_x(\mathcal{A}) \coloneqq \{u_{x}(A;\mathcal{A}) \colon A \in \mathcal{A}\}\textrm{, where }\]
\vspace{-0.05cm}
\[\hspace{0.1cm}u_x(A;\mathcal{A}) \coloneqq \begin{cases} A \cup \{x\} & A \cup \{x\} \not \in \mathcal{A} \\ A & A \cup \{x\} \in \mathcal{A} \end{cases}\textrm{.}\]

\vspace{0.3cm}

\noindent For each $S \in u_x(\mathcal{A})$, we also consider $u_x^{-1}(S;\mathcal{A}) \in \mathcal{A}$ such that $u_x(u_x^{-1}(S;\mathcal{A}); \mathcal{A}) = S$. In \cite{10}, two pertinent properties of up-compression were proved in the context of bounding average set size of a union-closed family. The first property is that $(u_p \circ \cdots \circ u_1)(\mathcal{A})$ is union-closed for any $p \in [n]$. Note that, since $u_x(\mathcal{A})$ has both the same size and same universe as of $\mathcal{A}$, we have by induction that $(u_p \circ \cdots \circ u_1)(\mathcal{A})$ also has the same size and universe as of $\mathcal{A}$. Now, call a family $\mathcal{F}$ with universe $[n]$ \textit{upward-closed} if it contains every superset of each of its members, i.e. if $X \in \mathcal{F}$ and $X \subseteq Y \subseteq [n]$ together imply that $Y \in \mathcal{F}$. The second property from \cite{10} applicable to our study is that the iteratively compressed family $(u_n \circ \cdots \circ u_1)(\mathcal{A})$ is upward-closed. We define $T_1$ to be this iterated compression:

\vspace{-0.05cm}

\[T_1 \coloneqq u_n \circ \cdots \circ u_1\textrm{.}\]

\medskip
\smallskip

\noindent We now demonstrate an additional property of up-compression, essential for proving the main theorem.

\medskip
\smallskip

\noindent \textbf{Lemma 2.1.1.} \textit{If} $\mathcal{A}$ \textit{is a union-closed family with universe} $[n]$\textit{, then for all} $x \in [n]$\textit{,} \[|\textsf{C}(u_x(\mathcal{A}),k)| \leq |\textsf{C}(\mathcal{A},k)|\textrm{.}\]

\begin{proof}
As up-compression bijectively maps members of $\mathcal{A}$ to members of $u_x(\mathcal{A})$, it suffices to show that if $\{C_1, \cdots, C_k\}$ is a $k$-chain in $u_x(\mathcal{A})$ such that $C_k \subseteq \cdots \subseteq C_1$, then $\{u_x^{-1}(C_1; \mathcal{A}), \cdots, u_x^{-1}(C_k;\mathcal{A})\}$ is a $k$-chain in $\mathcal{A}$ such that $u_x^{-1}(C_k;\mathcal{A}) \subseteq \cdots \subseteq u_x^{-1}(C_1;\mathcal{A})$. For any $i \in [k-1]$, we have as a result of $C_{i+1} \subseteq C_i$ that both $C_{i+1} \setminus \{x\} \subseteq C_i$ and $C_{i+1} \setminus \{x\} \subseteq C_i \setminus \{x\}$. Also, from the definition of up-compression, $u_x^{-1}(C_i; \mathcal{A}) \in \{C_i, C_i \setminus \{x\}\}$ and $u_x^{-1}(C_{i+1}; \mathcal{A}) \in \{C_{i+1}, C_{i+1} \setminus \{x\}\}$. Thus, if $u_x^{-1}(C_{i+1}; \mathcal{A}) = C_{i+1} \setminus \{x\}$ or $u_x^{-1}(C_{i+1}; \mathcal{A}) = C_{i+1} \ \land \ u_x^{-1}(C_i; \mathcal{A}) = C_i$, then $u_x^{-1}(C_{i+1}; \mathcal{A}) \subseteq u_x^{-1}(C_i; \mathcal{A})$. It remains to consider the possibility that $u_x^{-1}(C_i; \mathcal{A}) = C_i \setminus \{x\} \ \land \ u_x^{-1}(C_{i+1}; \mathcal{A}) = C_{i+1}$. In this case, if $x \in C_{i+1}$, then $x \in C_i$ and $u_x^{-1}(C_i; \mathcal{A}) \ \cup \ u_x^{-1}(C_{i+1}; \mathcal{A}) = (C_i \setminus \{x\}) \ \cup \ C_{i+1} = (C_i \setminus \{x\}) \ \cup \ (C_{i+1} \cup \{x\}) = C_i \cup C_{i+1} \cup \{x\} = C_i \cup \{x\} = C_i$. Because $\mathcal{A}$ is union-closed, it then follows that $C_i \in \mathcal{A}$. However, again via the definition of up-compression, $(x \in C_i \ \land \ C_i \in \mathcal{A}) \ \land \ (C_i \setminus \{x\} = u_x^{-1}(C_i; \mathcal{A}) \in \mathcal{A})$ implies that $u_x(u_x^{-1}(C_i; \mathcal{A}); \mathcal{A}) = C_i \setminus \{x\}$. Since $x \in C_i$, this contradicts that $u_x(u_x^{-1}(C_i; \mathcal{A}); \mathcal{A}) = C_i$. Therefore, it must be that $x \not \in C_{i+1}$, implying that $u_x^{-1}(C_{i+1}; \mathcal{A}) = C_{i+1} = C_{i+1} \setminus \{x\} \subseteq C_i \setminus \{x\} = u_x^{-1}(C_i; \mathcal{A})$.
\end{proof}

\medskip

\noindent Recalling that $\mathcal{A}_0$ is union-closed and that up-compression preserves the union-closed property, we inductively apply Lemma 2.1.1 in order to obtain that $|\textsf{C}(T_1(\mathcal{A}_0),k)| \leq |\textsf{C}(\mathcal{A}_0,k)|$.

\subsection*{2.2. Applying $T_2$ to $T_1(\mathcal{A}_0)$}

\noindent Let $\mathcal{A}$ again be a union-closed family with universe $[n]$. We denote the families consisting, respectively, of the smallest members of $\mathcal{A}$ and the largest members of the complement $2^{[n]} \setminus \mathcal{A}$ by the following:

\[\hspace{-0.5cm}\mathcal{S}_1(\mathcal{A}) \coloneqq \ \Bigr \{A \in \mathcal{A} \colon |A| = \min_{X \in \mathcal{A}}|X|\Bigr \}\textrm{,}\] 

\[\hspace{0.78cm}\mathcal{S}_2(\mathcal{A})\coloneqq \ \Bigr\{A \in 2^{[n]} \setminus \mathcal{A} \colon |A|=\max_{Y \in 2^{[n]} \setminus \mathcal{A}}|Y| \Bigr\}\textrm{.}\]

\smallskip

\noindent We set $\alpha_{\mathcal{A}}$ and $\beta_{\mathcal{A}}$ equal to the respective member sets lexicographically first in $\mathcal{S}_1(\mathcal{A})$ and $\mathcal{S}_2(\mathcal{A})$:

\[\hspace{0.2675cm} \alpha_\mathcal{A} = X \colon \ X \in \mathcal{S}_1(\mathcal{A}) \land (X' \in \mathcal{S}_1(\mathcal{A}) \setminus \{X\} \implies \min\{(X \cup X') \setminus (X \cap X')\} \in X)\textrm{,}\]

\medskip

\noindent \[\beta_\mathcal{A} = Y \colon \ Y \in \mathcal{S}_2(\mathcal{A}) \land (Y' \in \mathcal{S}_2(\mathcal{A}) \setminus \{Y\} \implies \min\{(Y \cup Y') \setminus (Y \cap Y')\} \in Y)\textrm{.}\]

\medskip

\noindent Next, we define an operator $V$ on $\mathcal{A}$ as follows:

\[V(\mathcal{A})= \begin{cases}(\mathcal{A} \setminus \{\alpha_{\mathcal{A}}\}) \cup \{\beta_{\mathcal{A}}\} & \ \textrm{if } |\alpha_\mathcal{A}| < |\beta_{\mathcal{A}}| \\ \mathcal{A} & \textrm{ if }|\alpha_\mathcal{A}| \geq |\beta_{\mathcal{A}}|\end{cases}\textrm{,}\]

\medskip

\noindent noting that $V(\mathcal{A})=\mathcal{A}$ if and only if $\mathcal{A}$ is top-aligned. Because $V(\mathcal{A}) \subseteq 2^{[n]}$ and $[n] \in V(\mathcal{A})$, we have that $U(V(\mathcal{A}))=U(\mathcal{A})=[n]$.

\bigskip

\noindent \textbf{Lemma 2.2.1.} \textit{If} $\mathcal{A}$ \textit{is upward-closed, then} $V(\mathcal{A})$ \textit{is also upward-closed.}

\begin{proof}
Assume $\mathcal{A}$ to be upward-closed. If $|\alpha_\mathcal{A}| \geq |\beta_{\mathcal{A}}|$, then $V(\mathcal{A})=\mathcal{A}$ and we are done. Thus, we let $|\alpha_\mathcal{A}| < |\beta_{\mathcal{A}}|$, as depicted in Figure 2.2.1, and consider any $X \in V(\mathcal{A})$. If $X \neq \beta_{\mathcal{A}}$, then we have by $V(\mathcal{A})\setminus \{\beta_{\mathcal{A}}\} = \mathcal{A} \setminus \{\alpha_{\mathcal{A}}\}$ that $X \in \mathcal{A}$, and it follows from $\mathcal{A}$ being upward-closed that $Y \in \mathcal{A}$ for any $Y \in 2^{[n]} \setminus\{X\}$ such that $X \subseteq Y$. Then, because $X \in \mathcal{A}$ implies that $|X| \geq |\alpha_{\mathcal{A}}|$, we have by $X \subseteq Y$ that $Y \neq \alpha_{\mathcal{A}}$, in turn implying that $Y \in \mathcal{A} \setminus \{\alpha_{\mathcal{A}}\} = V(\mathcal{A}) \setminus \{\beta_{\mathcal{A}}\} \subseteq V(\mathcal{A})$. Next, we assume that $X = \beta_{\mathcal{A}}$, and again consider any $Y \in 2^{[n]} \setminus\{X\}$ such that $X \subseteq Y$. In this case, $|\alpha_{\mathcal{A}}| < |\beta_{\mathcal{A}}| < |Y|$ implies that $Y \neq \alpha_{\mathcal{A}}$. Also, because $\beta_{\mathcal{A}}$ is of greatest size in $2^{[n]} \setminus \mathcal{A}$, we have that $Y \in \mathcal{A}$. It again follows that $Y \in \mathcal{A} \setminus \{\alpha_{\mathcal{A}}\} = V(\mathcal{A}) \setminus \{\beta_{\mathcal{A}}\} \subseteq V(\mathcal{A})$.
\end{proof}

\begin{figure}[h]

\caption*{\textbf{Figure 2.2.1:} If $|\alpha_{\mathcal{A}}| < |\beta_{\mathcal{A}}|$, then $V$ removes $\alpha_{\mathcal{A}}$ from $\mathcal{A}$ and adds $\beta_{\mathcal{A}}$ to the result, a transformation which preserves the upward-closed property.}

\vspace{-0.45cm}

\begin{center}

\begin{tikzpicture}[scale=1]

\def\x{8.5}
\def\y{0.5}

\node (A) at (2.35+\y,0.3) {$\mathcal{A}$};
\node (B) at (10.95,0.3) {$V(\mathcal{A})$};
\node[text=black] at (6.947,-1.2925) {$V$};

\fill[fill=gray!40, fill opacity=0.5] (\y,0) -- (\y,-2.25) -- (5+\y,0);
\fill[fill=gray!40, fill opacity=0.5] (\y,0) -- (3+\y,-2.5) -- (5+\y,0);
\fill[fill=gray!40, fill opacity=0.5] (\y,0) -- (5+\y,-1.2) -- (5+\y,0);
  
\node (1) at (5.825,-1.55) {};
\node (2) at (8.075,-1.55) {};
\draw[decoration={markings,mark=at position 1 with
{\arrow[scale=2,>=to]{>}}},postaction={decorate}] (1) -- (2);
  
\fill[fill=gray!40, fill opacity=0.5] (\y,0) -- (1.5+\y,-3.25) -- (5+\y,0) -- (\y,0) -- cycle;
\fill[fill=gray!40, fill opacity=0.5] (\x,0) -- (\x,-2.25) -- (\x+5,0);
\fill[fill=gray!40, fill opacity=0.25] (\x,0) -- (\x+1.35,-2.95) -- (\x+5,0);
\fill[fill=gray!40, fill opacity=0.25] (\x,0) -- (\x+1.5,-2.95) -- (\x+5,0);
\fill[fill=gray!40, fill opacity=0.25] (\x,0) -- (\x+1.65,-2.95) -- (\x+5,0);
\fill[fill=gray!40, fill opacity=0.5] (\x,0) -- (\x+3,-2.5) -- (\x+5,0);
\fill[fill=gray!40, fill opacity=0.5] (\x,0) -- (\x+5,-1.2) -- (\x+5,0);
\fill[fill=gray!40, fill opacity=0.5] (\x,0) -- (12.7125,-1.15) -- (\x+5,0) -- (\x,0) -- cycle;
    
\node[text=black] at (1.875+\y,-3.25) {$\alpha_{\mathcal{A}}$};
\node[scale = 0.375, circle, fill=black] at (1.5+\y,-3.25) {};
\node[text=black] at (12.7,-1.42) {$\beta_{\mathcal{A}}$};
\node[scale = 0.375, circle, fill=black] at (12.7125,-1.15) {};

\end{tikzpicture}

\end{center}

\end{figure}
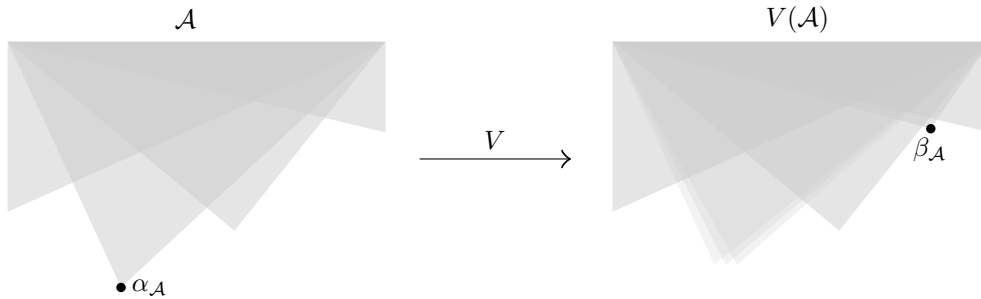

\vspace{-0.25cm}

\noindent \textbf{Lemma 2.2.2.} \textit{If} $\mathcal{A}$ \textit{is upward-closed and non-top-aligned with size $m$}\textit{, then} $|\textsf{C}(V(\mathcal{A}),k)| \leq |\textsf{C}(\mathcal{A},k)|$\textit{. If, in addition,} $c(k,m,n) > 0$\textit{, then} $|\textsf{C}(V(\mathcal{A}),k)| < |\textsf{C}(\mathcal{A},k)|$\textit{.}

\vspace{0.0125cm}

\begin{proof}
Assume $\mathcal{A}$ to be upward closed and non-top-aligned ($|\alpha_{\mathcal{A}}| < |\beta_{\mathcal{A}}|$) with $|\mathcal{A}|=m$. We first define 

\vspace{-0.0375cm}

\[\hspace{-0.25cm}\textsf{C}_{\alpha}(\mathcal{A},k) \coloneqq \{\mathcal{C} \in \textsf{C}(\mathcal{A},k) \colon \alpha_{\mathcal{A}} \in \mathcal{C})\}\textrm{ and}\]
\[\hspace{0.225cm}\textsf{C}_{\beta}(V(\mathcal{A}),k) \coloneqq \{\mathcal{C} \in \textsf{C}(V(\mathcal{A}),k) \colon \beta_{\mathcal{A}} \in \mathcal{C})\}\textrm{,}\] 

\medskip

\noindent observing that $\textsf{C}(\mathcal{A},k) \setminus \textsf{C}_{\alpha}(\mathcal{A},k) = \textsf{C}(V(\mathcal{A}),k) \setminus \textsf{C}_{\beta}(V(\mathcal{A}),k)$. Based on this identity, if we show that $|\textsf{C}_{\beta}(V(\mathcal{A}),k)| \leq |\textsf{C}_{\alpha}(\mathcal{A},k)|$, then we will have proved that $|\textsf{C}(V(\mathcal{A}),k)| \leq |\textsf{C}(\mathcal{A},k)|$. To establish that $|\textsf{C}_{\beta}(V(\mathcal{A}),k)| \leq |\textsf{C}_{\alpha}(\mathcal{A},k)|$, we first note that $\mathcal{A}$ being upward-closed implies that there exists $Y \in \mathcal{A}$ lexicographically first such that $\alpha_{\mathcal{A}} \subseteq Y$ and $|Y| = |\beta_{\mathcal{A}}|$. Also, the upward-closed property of $\mathcal{A}$, together with the fact that $\beta_{\mathcal{A}} \not \in \mathcal{A}$, implies that $\mathcal{A}$ does not contain any proper subset of $\beta_{\mathcal{A}}$, in turn implying that neither does $V(\mathcal{A})$. It follows that $C_k = \beta_{\mathcal{A}}$ for any $\mathcal{C}=\{C_1, \cdots, C_k\} \in \textsf{C}_{\beta}(V(\mathcal{A}),k)$ such that $C_k \subseteq \cdots \subseteq C_1$. Since $|Y|=|\beta_{\mathcal{A}}|$, and all proper supersets of $\beta_{\mathcal{A}}$ and proper supersets of $Y$ belong to $\mathcal{A} \cap V(\mathcal{A})$, each $\mathcal{C}=\{C_1, \cdots, C_{k-1}, \beta_{\mathcal{A}}\} \in \textsf{C}_{\beta}(V(\mathcal{A}),k)$ must have (as a result of the natural bijection from $[n] \setminus \beta_{\mathcal{A}}$ to $[n] \setminus Y$) a corresponding $k$-chain $\mathcal{C}' = \{C'_1, \cdots, C'_{k-1}, Y\} \in \textsf{C}(\mathcal{A},k)$ such that $Y \subseteq C'_{k-1} \subseteq \cdots \subseteq C'_1$. Then, because $\alpha_{\mathcal{A}}$ is a proper subset of $Y$ in $\mathcal{A}$, we have that for every such $\mathcal{C}'$, there is a unique $k$-chain $\mathcal{C}''= (\mathcal{C}'\setminus \{Y\}) \cup \{\alpha_{\mathcal{A}}\} = \{C_1', \cdots, C'_{k-1}, \alpha_{\mathcal{A}}\} \in \textsf{C}_{\alpha}(\mathcal{A},k)$ such that $\alpha_{\mathcal{A}} \subseteq C'_{k-1} \subseteq \cdots \subseteq C'_1$. As the map $\Psi \colon \textsf{C}_{\beta}(V(\mathcal{A}),k) \to \textsf{C}_{\alpha}(\mathcal{A},k)$ given by $\Psi(\mathcal{C}) = \mathcal{C}''$ is injective, we obtain that $|\textsf{C}_{\beta}(V(\mathcal{A}),k)| \leq |\textsf{C}_{\alpha}(\mathcal{A},k)|$, which, by the aforementioned identity, implies that $|\textsf{C}(V(\mathcal{A}),k)| \leq |\textsf{C}(\mathcal{A},k)|$. Finally, to establish that this inequality is strict whenever $c(k,m,n) > 0$, we show that the collection $\textsf{C}_\alpha(\mathcal{A},k) \setminus \Psi(\textsf{C}_{\beta}(V(\mathcal{A}),k))$ under this additional assumption is nonempty. Since $\alpha_{\mathcal{A}}$ is of minimum size in $\mathcal{A}$, we have that $m \leq \sum_{i=0}^{n-|\alpha_{\mathcal{A}}|} \binom{n}{i}$. From this, it follows that $k \leq n - |\alpha_{\mathcal{A}}|+1$. (Otherwise, no top-aligned family with such size and universe would have a $k$-chain, contradicting that $c(k,m,n)>0$.) Because $\mathcal{A}$ is upward-closed, it then follows that there exist a pair of nonnegative integers $(s,t)$ and chain $\mathcal{C}^* = \{C^{\supseteq Y}_1, \cdots, C^{\supseteq Y}_s, Y, C^{\subseteq Y}_1, \cdots, C^{\subseteq Y}_t, \alpha_{\mathcal{A}}\} \in \textsf{C}_{\alpha}(\mathcal{A},k)$ such that $\alpha_{\mathcal{A}} \subseteq C^{\subseteq Y}_t \subseteq \cdots \subseteq C^{\subseteq Y}_1 \subseteq Y \subseteq C^{\supseteq Y}_s \subseteq \cdots \subseteq C^{\supseteq Y}_1$ and $s+t+2=k$. Since $Y \in \mathcal{C}^*$ and $Y \not \in \mathcal{C}''$ for all $\mathcal{C}'' \in \Psi(\textsf{C}_{\beta}(V(\mathcal{A}),k))$, it must be that $\mathcal{C}^* \in \textsf{C}_\alpha(\mathcal{A},k) \setminus \Psi(\textsf{C}_{\beta}(V(\mathcal{A}),k))$.
\end{proof}

\medskip

\noindent Let $V^0(T_1(\mathcal{A}_0)) = T_1(\mathcal{A}_0)$, and for any nonnegative integer $i$, let $V^{i+1}(T_1(\mathcal{A}_0)) = V(V^i(T_1(\mathcal{A}_0)))$. Also, let $z$ be the smallest integer such that $V^{z}(T_1(\mathcal{A}_0)) = V^{z+1}(T_1(\mathcal{A}_0))$, i.e. the number of consecutive applications of $V$ required to transform $T_1(\mathcal{A}_0)$ into a top-aligned family. We define $T_2$ to be this iterated application of $V$:

\vspace{-0.4cm}

\[T_2 \coloneqq V^z\textrm{.}\]

\vspace{0.15cm}

\noindent By Lemma 2.2.1, if $V^i(T_1(\mathcal{A}_0))$ is upward-closed, then $V^{i+1}(T_1(\mathcal{A}_0))$ is also upward-closed. Since $V^0(T_1(\mathcal{A}_0))$ is upward-closed, it follows by induction that $V^i(T_1(\mathcal{A}_0))$ is upward-closed for all $i \geq 0$. When $i < z$, we then have by Lemma 2.2.2 that $|\textsf{C}(V^{i+1}(T_1(\mathcal{A}_0)),k)| \leq |\textsf{C}(V^i(T_1(\mathcal{A}_0)),k)|$. This implies, again by induction, that $|\textsf{C}(V^z(T_1(\mathcal{A}_0)),k)| \leq |\textsf{C}(V^0(T_1(\mathcal{A}_0)),k)| = |\textsf{C}(T_1(\mathcal{A}_0),k)|$. Since $|\textsf{C}(\mathcal{A}_{\min},k)| = |\textsf{C}(T_2(T_1(\mathcal{A}_0)),k)| = |\textsf{C}(V^z(T_1(\mathcal{A}_0)),k)|$ and $|\textsf{C}(T_1(\mathcal{A}_0),k)| \leq |\textsf{C}(\mathcal{A}_0,k)|$, it finally follows that $|\textsf{C}(\mathcal{A}_{\min},k)| \leq |\textsf{C}(\mathcal{A}_0,k)|$, which establishes the first part of our main theorem.

\subsection*{2.3. Uniqueness}

Having proved that $c(k,m,n)$ is achieved by top-aligned families with universe $[n]$ and size $m$, we now transition to proving the second part of the theorem, namely that top-aligned families are in fact the only minimizers whenever the minimum is nonzero and $m > n$. We first show that uniqueness does not hold anywhere in the regime $c(k,m,n)>0 \ \land \ 2 \leq m \leq n$. (Here, $k=2$, as $c(k,m,n)>0$ implies that any top-aligned family with universe $[n]$ and size $m$ must have at least one $k$-chain, yet   $m \leq n$ implies that top-aligned families only contain member sets from the top two layers of the power set.)

\medskip 

\noindent \textbf{Proposition 2.3.1.} \textit{For all integers} $m$ \textit{and} $n$ \textit{such that} $2 \leq m \leq n$\textit{, there exists a non-top-aligned union-closed family} $\tilde{\mathcal{A}}_{m,n}$ \textit{with universe} $[n]$ \textit{and size} $m$ \textit{such that} $|\textsf{C}(\tilde{\mathcal{A}}_{m,n},2)| = c(2,m,n) > 0$\textit{.} 

\begin{proof}Set $[0] = \emptyset$. For all integers $m$ and $n$ such that $2 \leq m \leq n$, we let \[\tilde{\mathcal{A}}_{m,n} = \Bigr \{[n],[m-2]\Bigr\} \ \cup \ \Bigr\{[n] \setminus \{x\} \colon x \in [m-2]\Bigr\}\textrm{.}\] The family $\tilde{\mathcal{A}}_{m,n}$ is union-closed with universe $[n]$ and size $m$, yet not top-aligned because $[m-2]$ is a member set of $\tilde{\mathcal{A}}_{m,n}$ smaller than the set $[n] \setminus \{m\}$ from $2^{[n]} \setminus \tilde{\mathcal{A}}_{m,n}$. A key property of $\tilde{\mathcal{A}}_{m,n}$ is that $X \subseteq Y$ implies that $Y=[n]$ for any distinct $X$ and $Y$ therein. Since $[n] \in \tilde{\mathcal{A}}_{m,n}$, it follows that $|\textsf{C}(\tilde{\mathcal{A}}_{m,n},2)| = m-1$. Finally, we note that any top-aligned family with universe $[n]$ and size $m$ such that $2 \leq m \leq n$ also has $m-1$ pairs of comparable member sets. Since top-aligned families achieve the minimum $c(2,m,n)$, this implies that $|\textsf{C}(\tilde{\mathcal{A}}_{m,n},2)| = c(2,m,n) > 0$.\end{proof}

\medskip

\noindent Proceeding with the main proof, we now assume that $m>n$ and $c(k,m,n)>0$. By Lemma 2.2.2, $\mathcal{A}_{\min}$ has strictly fewer $k$-chains than $\mathcal{A}_0$ if $T_1(\mathcal{A}_0)$ is not yet top-aligned. Hence, we need only consider the case where $T_1(\mathcal{A}_0)$ is already top-aligned (rendering $T_2$ to act as the identity operator $V^0$).

\medskip

\noindent For any union-closed family $\mathcal{A}$, let $u_0(\mathcal{A})=\mathcal{A}$. We define the following truncation of $T_1$: \[T_1^- \coloneqq u_{r-1} \circ \cdots \circ u_0\textrm{,}\] 
\noindent where $r$ is greatest in $[n]$ such that $(u_{r-1} \circ \cdots \circ u_0)(\mathcal{A}_0) \neq (u_r \circ \cdots \circ u_0)(\mathcal{A}_0)$. 

\medskip

\noindent We first suppose that $m=n+1$, which implies that $T_1(\mathcal{A}_0)$, being top-aligned, is the union of the top two layers of the power set. In this case, $c(k,m,n)>0$ implies that $k=2$ (as in Proposition 2.3.1), and we have that $c(2,m,n)=|\textsf{C}(T_1(\mathcal{A}_0),2)|=m-1$. Further, from the definition of up-compression, $[n] \setminus \{r\} \in T_1(\mathcal{A}_0)$ implies that $u_r^{-1}([n] \setminus \{r\}; T_1^-(\mathcal{A}_0)) = [n] \setminus \{r\} \in T_1^-(\mathcal{A}_0)$. However, because $T_1(\mathcal{A}_0)=u_r(T_1^-(\mathcal{A}_0)) \neq T_1^-(\mathcal{A}_0)$, there also exists $X \in T_1^-(\mathcal{A}_0)$ such that $|X| = n-2$ and $r \not \in X$, and it follows that $X$ is a proper subset of $[n] \setminus \{r\}$ in $T_1^-(\mathcal{A}_0)$. Since $|T_1^-(\mathcal{A}_0)|=|T_1(\mathcal{A}_0)|$, and $[n] \in T_1^-(\mathcal{A}_0)$ with $A \subseteq [n]$ for all $A \in T_1^-(\mathcal{A}_0) \setminus \{[n]\}$, it then follows that $|\textsf{C}(T_1^-(\mathcal{A}_0),2)|>m-1$. As $|\textsf{C}(T_1^-(\mathcal{A}_0),2)| \leq |\textsf{C}(\mathcal{A}_0,2)|$ and $c(2,m,n)=m-1$, this then implies that $|\textsf{C}(\mathcal{A}_0,2)| > c(2,m,n)$.

\medskip

\noindent Next we assume that $m>n+1$. Let $q$ be the greatest integer such that $T_1(\mathcal{A}_0) \cap \{X \subseteq [n] \colon |X| \leq q\} = \emptyset$. If there exists some $A \in 2^{[n]} \setminus T_1^-(\mathcal{A}_0)$ with $|A| \geq q+2$ and $r \not \in A$, then $A \not \in u_r(T_1^-(\mathcal{A})) = T_1(\mathcal{A}_0)$, which contradicts that $T_1(\mathcal{A}_0)$, again being top-aligned, must contain all possible sets of size at least $q+2$. Thus, for any $A \in 2^{[n]}$, we have that if $|A|\geq q+2$ and $r \not \in A$, then $A \in T_1^-(\mathcal{A}_0)$. This in turn implies that any $A \in 2^{[n]}$ with $|A| \geq q+3$ must belong to $T_1^-(\mathcal{A}_0)$. (If not, then there exists some $A \in 2^{[n]} \setminus T_1^-(\mathcal{A}_0)$ such that $|A|\geq q+3$ and $r \in A$, and it follows from $A \setminus \{r\} \in T_1^-(\mathcal{A}_0)$ that $u_r(A \setminus \{r\}; T_1^-(\mathcal{A}_0)) = A \in T_1(\mathcal{A}_0)$. This implies that $A \setminus \{r\} \not \in T_1(\mathcal{A}_0)$, again contradicting that $T_1(\mathcal{A}_0)$ must contain all possible sets of size at least $q+2$.) As $T_1^-(\mathcal{A}_0) \neq T_1(\mathcal{A}_0)$ and $|u_r^{-1}(A;T_1^-(\mathcal{A}_0))| \geq |A|-1$ for all $A \in T_1(\mathcal{A}_0)$, we then have that there is some $Y \in T_1^-(\mathcal{A}_0)$ with $q \leq |Y| \leq q+1$ such that $r \not \in Y$ and $Y \cup \{r\} \not \in T_1^-(\mathcal{A}_0)$. Since $q \leq n-3$ (implied by the assumption that $m> n+1$), there also exists some $Z \in T_1^-(\mathcal{A}_0)$ of size $q+2$ such that $r \not \in Z$. Now, since $m \leq \sum_{i=0}^{n-q-1} \binom{n}{i}$, we have that $k \leq n-q$. (Otherwise, any top-aligned family with such size and universe would be without a $k$-chain, contradicting that $c(k,m,n)>0$.) This then implies that there exists a $k$-chain $\{\hat{C}_1, \cdots, \hat{C}_{k-2}, Z, Y\} \in \textsf{C}(T_1^-(\mathcal{A}_0),k)$ such that $Y \subseteq Z \subseteq \hat{C}_{k-2} \subseteq \cdots \subseteq \hat{C}_1$. Because $T_1^-(\mathcal{A}_0)$ contains $Z$ as well as $Z \cup \{r\}$, we have that $u_r(Z; T_1^-(\mathcal{A}_0))=Z$. Therefore, $r \not \in u_r(Z; T_1^-(\mathcal{A}_0))$, which, together with $u_r(Y;T_1^-(\mathcal{A}_0)) = Y \cup \{r\}$, implies that $u_r(Y;T_1^-(\mathcal{A}_0)) \not \subseteq u_r(Z; T_1^-(\mathcal{A}_0))$. Recall that, in proving Lemma 2.1.1, we showed that for any union-closed family $\mathcal{A}$ and $x \in [n]$, if $\{C_1, \cdots, C_k\}$ is a $k$-chain in $u_x(\mathcal{A})$ such that $C_k \subseteq \cdots \subseteq C_1$, then $\{u_x^{-1}(C_1; \mathcal{A}), \cdots, u_x^{-1}(C_k;\mathcal{A})\}$ is a $k$-chain in $\mathcal{A}$ such that $u_x^{-1}(C_k;\mathcal{A}) \subseteq \cdots \subseteq u_x^{-1}(C_1;\mathcal{A})$. Since $\{\hat{C}_1, \cdots, \hat{C}_{k-2}, Z, Y\}$ is a $k$-chain in $T_1^{-1}(\mathcal{A}_0)$ such that $Y \subseteq Z \subseteq \hat{C}_{k-2} \subseteq \cdots \subseteq \hat{C}_1$, yet $\{u_r(\hat{C}_1; T_1^{-1}(\mathcal{A}_0)), \cdots, u_r(\hat{C}_{k-2}; T_1^{-1}(\mathcal{A}_0)), u_r(Z; T_1^{-1}(\mathcal{A}_0)), u_r(Y; T_1^{-1}(\mathcal{A}_0))\}$ is not a corresponding $k$-chain in $u_r(T_1^{-1}(\mathcal{A}_0))$  (in particular, $u_r(Y;T_1^-(\mathcal{A}_0)) \not \subseteq u_r(Z;T_1^-(\mathcal{A}_0))$), it follows that $|\textsf{C}(T_1^-(\mathcal{A}_0),k)| > |\textsf{C}(T_1(\mathcal{A}_0,k))|$. As $|\textsf{C}(T_1^-(\mathcal{A}_0),k)| \leq |\textsf{C}(\mathcal{A}_0,k)|$ and $|\textsf{C}(T_1(\mathcal{A}_0),k)|=c(k,m,n)$, we finally obtain that $|\textsf{C}(\mathcal{A}_0,k)| > c(k,m,n)$. 

\medskip

\noindent Hence, top-aligned families are uniquely extremal when $m>n$ and $c(k,m,n)>0$, which concludes the proof of the main theorem.

\section*{3. A corollary of the theorem}

A family of sets $\mathcal{F}$ with universe $[n]$ is \textit{separating} if, for all $x$ and $y$ in $[n]$, there exists a member set $F$ in $\mathcal{F}$ such that $(x \in F) \land (y \not \in F) \lor (x \not \in F) \land (y \in F)$.

\medskip

\noindent \textbf{Corollary.} \textit{For any positive integers} $k$ \textit{and} $n$ \textit{such that} $2 \leq k \leq n+1$\textit{, the number of} $k$\textit{-chains in a separating union-closed family with universe} $[n]$ \textit{and size} $m$ \textit{is minimized whenever the family is top-aligned. Further, if the minimum is nonzero, then there are no other minimizing families.}

\begin{proof} We first note that (by Lemma 2 of \cite{7}) the size of any separating union-closed family is greater than or equal to the size of its universe. Additionally, we observe that any top-aligned family with size at least that of its universe is guaranteed to be separating. It follows by the main theorem that for $2 \leq k \leq n+1$, the minimum number of $k$-chains in a separating union-closed family with universe $[n]$ and size $m$ is indeed equal to $c(k,m,n)$, achieved whenever the family is top-aligned. The main theorem also then implies that there are no other extremal (separating) families whenever $c(k,m,n) > 0$ and $m>n$. When $m=n$, we again have by $c(k,m,n)>0$ that $k = 2$. Thus, it remains only to prove that any separating union-closed family $\mathcal{A}_{\min}$ with universe $[n]$ and size $m$ such that $m=n$ and $|\textsf{C}(\mathcal{A}_{\min},2)|=c(2,m,n)=m-1>0$ must be top-aligned. To this end, we observe that for every $x \in [n]$, there is at most one $A \in \mathcal{A}_{\min}$ such that $x \not \in A$. (Otherwise, there exist distinct sets $A_1,A_2 \in \mathcal{A}_{\min}$ and an element $x \in [n]$ such that $x \not \in A_1$ and $x \not \in A_2$. This implies that $A_1 \cup A_2 \neq [n]$, yet $A_1 \cup A_2 \in \mathcal{A}_{\min}$ by the union-closed property. Since $[n] \in \mathcal{A}_{\min}$, and $A \subseteq [n]$ for every $A \in \mathcal{A}_{\min} \setminus \{[n]\}$, it then follows from $(A_1 \subseteq A_1 \cup A_2 \ \land \ A_1 \neq A_1 \cup A_2) \ \lor \ (A_2 \subseteq A_1 \cup A_2 \ \land \ A_2 \neq A_1 \cup A_2)$ that $|\textsf{C}(\mathcal{A}_{\min},2)| > m-1$, a contradiction.) As a result, $\sum_{x \in [n]} |\{A \in \mathcal{A}_{\min} \colon x \in A\}| \geq mn - n = n(m-1)$. In addition, we have by double counting that $\sum_{x \in [n]}|\{A \in \mathcal{A}_{\min} \colon x \in A\}| = \sum_{A \in \mathcal{A}_{\min}}|A|$. Hence, $\sum_{A \in \mathcal{A}_{\min}}|A| \ \geq \ n(m-1) = m(n-1)$, which implies that the average set size of $\mathcal{A}_{\min}$ is at least $n-1$. It follows that either (\romannumeral 1) every member set of $\mathcal{A}_{\min} \setminus \{A,[n]\}$ has size $n-1$, where $A \in \mathcal{A}_{\min}$ such that $|A|=n-2$, or (\romannumeral 2) every member set of $\mathcal{A}_{\min} \setminus \{[n]\}$ has size $n-1$. If condition (\romannumeral 1) holds, then let $x_1$ and $x_2$ be the two distinct elements of $[n] \setminus A$. Because $A$ is a subset of both $A \cup \{x_1\}$ and $A \cup \{x_2\}$, if $A \cup \{x_1\} \in \mathcal{A}_{\min}$ or $A \cup \{x_2\} \in \mathcal{A}_{\min}$, then $|\textsf{C}(\mathcal{A}_{\min},k)| > m-1$, which contradicts that $|\textsf{C}(\mathcal{A}_{\min},k)|=c(2,m,n) = m-1$. Thus, $A \cup \{x_1\} \not \in \mathcal{A}_{\min}$ and $A \cup \{x_2\} \not \in \mathcal{A}_{\min}$, and it follows that every member set of $\mathcal{A}_{\min} \setminus \{A\}$ contains both $x_1$ and $x_2$, which then contradicts the assumption that $\mathcal{A}_{\min}$ is separating. Therefore, condition (\romannumeral 2) holds, which is equivalent to $\mathcal{A}_{\min}$ being top-aligned. \end{proof}


\vspace{0.05cm}

\begin{thebibliography}{12}

\bibitem{1}
J. Balogh and A.Z. Wagner, \textit{Kleitman's conjecture about families of given size minimizing the number of} $k$\textit{-chains}, Adv. Math. \textbf{330} (2018), 229--252.
\bibitem{2}
C. Bouchard, \textit{An upper bound for union-closed family size}, preprint (2025), arXiv:2511.10608.
\bibitem{3}
H. Bruhn and O. Schaudt, \textit{The journey of the union-closed sets conjecture}, Graphs Combin. \textbf{31} (2015), 2043--2074.
\bibitem{4}
S. Das, W. Gan, and B. Sudakov, \textit{Sperner's theorem and a problem of Erd\H os, Katona and Kleitman}, Combin. Probab. Comput. \textbf{24} (2015), 585--608.
\bibitem{5}
A.P. Dove, J.R. Griggs, R.J. Kang, and J.-S. Sereni, \textit{Supersaturation in the Boolean lattice}, Integers \textbf{14A} (2014), A4.
\bibitem{6}
P. Erd\H os, \textit{On a lemma of Littlewood and Offord}, Bull. Amer. Math. Soc. \textbf{51} (1945), 898--902.
\bibitem{7}
V. Falgas-Ravry, \textit{Minimal weight in union-closed families}, Electron. J. Combin. \textbf{18(1)} (2011), P95.
\bibitem{8}
J. Gilmer, \textit{A constant lower bound for the union-closed sets conjecture}, preprint (2022), arXiv:2211.09055.
\bibitem{9}
D. Kleitman, \textit{A conjecture of Erd\H os-Katona on commensurable pairs among subsets of an n-set (Theory of Graphs: Proc. Colloq. Tihany 1966)}, Academic Press (1968), 215--218.
\bibitem{10}
D. Reimer, \textit{An average set size theorem}, Combin. Probab. Comput. \textbf{12} (2003), 89--93.
\bibitem{11}
W. Samotij, \textit{Subsets of posets minimising the number of chains}, Trans. Amer. Math. Soc. \textbf{371} (2019), 7259--7274.
\bibitem{12}
E. Sperner, \textit{Ein Satz \"uber Untermengen einer endlichen Menge}, Math. Z. \textbf{27} (1928), 544--548.

\end{thebibliography}
\end{document}